\documentclass[final,3p,times,authoryear]{elsarticle}
\usepackage{amsmath,amssymb,amsthm,bm,mathtools}
\usepackage{microtype}
\usepackage{flafter}
\usepackage[colorlinks=true,allcolors=blue]{hyperref}
\makeatletter
\let\ps@pprintTitle\ps@plain
\makeatother
\journal{Statistics \& Probability Letters}
\newtheorem{theorem}{Theorem}[section]
\newtheorem{proposition}[theorem]{Proposition}
\newtheorem{lemma}[theorem]{Lemma}
\theoremstyle{definition}
\newtheorem{example}[theorem]{Example}
\theoremstyle{remark}

\newcommand{\R}{\mathbb R}
\newcommand{\E}{\mathbb E}
\newcommand{\Prb}{\mathbb P}
\newcommand{\M}{\mathcal M}
\newcommand{\supp}{\operatorname{supp}}
\newcommand{\conv}{\operatorname{conv}}

\newcommand{\diam}{\operatorname{diam}}
\newcommand{\argmin}{\operatorname*{arg\,min}}
\newcommand{\norm}[1]{\lVert#1\rVert}
\newcommand{\abs}[1]{\lvert#1\rvert}
\hypersetup{pdftitle={Unconstrained support points can leave a convex support}}
\newif\ifanonymous
\anonymousfalse
\ifanonymous\else
\hypersetup{pdfauthor={Zhaohui Li}}
\fi

\begin{document}
\begin{frontmatter}
\title{Unconstrained support points can leave a convex support}
\ifanonymous
\else
\author{Zhaohui Li}
\ead{lizh@amss.ac.cn}
\address{State Key Laboratory of Mathematical Sciences (SKLMS), Academy of Mathematics and Systems Science, Chinese Academy of Sciences, Beijing 100190, China}
\fi
\begin{abstract}
Support points provide a practical and effective way to approximate a probability distribution by a small set of equally weighted representatives through energy-distance minimization.
In many applications, the variables satisfy convex constraints that the representatives should also respect.
We investigate whether unconstrained support points automatically satisfy such constraints and show that this is not guaranteed.
We construct counterexamples using three planar observations, a density with full triangular support, and positive-definite matrix data admitting an indefinite optimal representative under the Frobenius norm.
We also identify cases in which containment holds, establish stability under distributional perturbations, and quantify the objective cost of imposing a convex-hull constraint.
For a fixed number of representatives and a positive population constraint gap, empirical configurations with vanishing global objective error eventually leave the sample convex hull almost surely.
\end{abstract}
\end{frontmatter}

\section{Introduction and setup}
Support points provide a practical way to represent a probability distribution by a small set of equally weighted points chosen through energy-distance minimization \citep{MakJoseph2018}.
Other representative-point methods include principal points, which minimize the expected squared distance to the nearest representative \citep{Flury1990}, and kernel herding \citep{ChenWellingSmola2010}.
The energy distance used by support points is proportional to squared maximum mean discrepancy (MMD) for a suitable kernel \citep{SejdinovicEtAl2013,MuandetEtAl2017}, connecting support points to kernel-based methods.
Projected support points address high-dimensional data reduction \citep{MakJoseph2017}, while SPlit, Data Twinning, MMD quantisation over observed candidates and kernel thinning construct representative subsets of the observations \citep{JosephVakayil2022,VakayilJoseph2022,TeymurEtAl2021,DwivediMackey2024}.

Many sampling problems involve convex constraints, including sampling from convex bodies \citep{LovaszVempala2006} and simulating normal distributions under linear restrictions \citep{Botev2017}.
When support points are used to represent such distributions, it is natural to require that they respect the same constraints.
We ask whether unconstrained support points automatically lie in the convex feasible region containing the observations.
A direct conjecture is that they lie in the sample convex hull.
Indeed, if $y_1,\ldots,y_N\in K$ and $K\subseteq\R^d$ is convex, then $\conv\{y_1,\ldots,y_N\}\subseteq K$.
Sample-hull containment is therefore a sufficient, potentially stronger condition for feasibility, since the hull may be strictly smaller than $K$.

We give counterexamples to this containment conjecture, starting with three planar observations and extending the construction to a density with full triangular support.
Leaving the sample hull alone does not imply violation of the original constraints.
However, the full-support example shows that representatives can leave the convex support itself, and a further example with positive-definite matrix data admits an indefinite optimal representative under the Frobenius norm.
We also study whether more observations restore containment and quantify the objective cost of a hull constraint.

In the original formulation, support points belong to a specified sample space $\mathcal X$ \citep{MakJoseph2018}.
Choosing the feasible region as $\mathcal X$ imposes the constraint explicitly.
Here we study whether it can be omitted by taking the optimization domain to be $\R^d$.

Let $d\ge1$ and $P$ be a Borel probability measure on $\R^d$ with finite first moment.
Fix an integer $m\ge1$ and write $\bm x=(x_1,\ldots,x_m)\in(\R^d)^m$.
Let $Y,Y'$ be independent with law $P$, and let $\delta_x$ denote the unit point mass at $x$.
Omitting the constant $-\E\norm{Y-Y'}$, the energy-distance objective for approximating $P$ by $m^{-1}\sum_i\delta_{x_i}$ is
\begin{equation}
J_P(\bm x)=\frac2m\sum_{i=1}^m\E\norm{x_i-Y}
-\frac1{m^2}\sum_{i,j=1}^m\norm{x_i-x_j}.
\label{eq:objective}
\end{equation}
Here $\norm{\cdot}$ is the Euclidean norm; see \citet{BaringhausFranz2004,SzekelyRizzo2013} for the related distance-based statistics.
For any nonempty set $D\subseteq\R^d$, write $D^m=D\times\cdots\times D$ for its $m$-fold Cartesian product and define
\[
\M_m(P)=\argmin_{\bm x\in(\R^d)^m}J_P(\bm x),\quad
V_m(P;D)=\inf_{\bm x\in D^m}J_P(\bm x),\quad
V_m(P)=V_m(P;\R^d).
\]
Write $\supp P$ for the topological support of $P$, $\conv(D)$ for the convex hull of $D$, and $\overline{\conv}(D)$ for its closure.
Thus $\M_m(P)$ is the set of all global minimizers, while $V_m(P;D)$ and $V_m(P)$ are the constrained and unconstrained infima, respectively.
For a nonempty compact set $C\subseteq\R^d$ containing $\supp P$, write
\begin{equation}
\Gamma_m(P,C)=V_m(P;C)-V_m(P)\ge0
\label{eq:price}
\end{equation}
for the increase in the minimum objective caused by the constraint.
A positive gap certifies that every unconstrained global minimizer has at least one point outside $C$; we use this criterion to establish the counterexamples below.

\section{Containment and counterexamples}
Containment holds for one representative or for one-dimensional data.

\begin{proposition}[Cases with automatic containment]
\label{prop:safe}
Let $C=\overline{\conv}(\supp P)$.
If $m=1$ or $d=1$, every global minimizer lies in $C^m$.
In dimension one, a tuple sorted as $x_1\le\cdots\le x_m$ is optimal if and only if
\begin{equation}
F(x_i-)\le \frac{2i-1}{2m}\le F(x_i),\qquad i=1,\ldots,m,
\label{eq:quantiles}
\end{equation}
where $F(t)=P(( -\infty,t])$ and $F(t-)=P(( -\infty,t))$.
\end{proposition}
\begin{proof}
For $m=1$, projection onto $C$ strictly decreases the expected distance at any exterior point.
For $d=1$, sorting does not change the objective, and sorted tuples satisfy $\sum_{i,j}\abs{x_i-x_j}=2\sum_i(2i-m-1)x_i$.
Thus \eqref{eq:objective} separates into convex scalar terms with subgradient intervals $\frac4m[F(x_i-)-u_i,F(x_i)-u_i]$, where $u_i=(2i-1)/(2m)$.
Their minimum sets are ordered, satisfy \eqref{eq:quantiles}, and lie in $C$; see also \citet{MakJoseph2018}.
\end{proof}

\begin{figure*}[t]
\centering
\includegraphics[width=\textwidth]{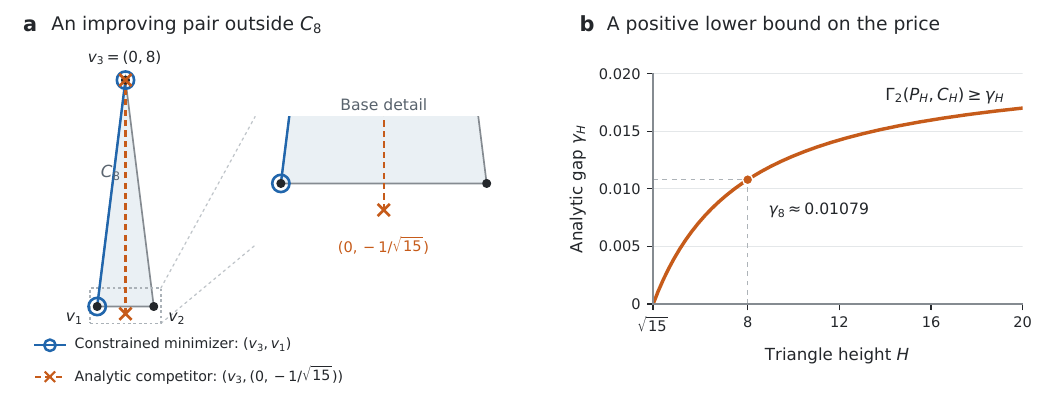}
\caption{Left: $C_8$, a hull-constrained optimal pair and the exterior competitor share the top vertex.
Right: their objective difference $\gamma_H$, a lower bound on $\Gamma_2(P_H,C_H)$.
The unconstrained minimum is not computed.}
\label{fig:counterexample}
\end{figure*}

The preceding proposition shows that failure requires at least two representatives in at least two dimensions.
We start with this case.
Figure~\ref{fig:counterexample} illustrates the construction.

\begin{example}[Three planar observations]
\label{ex:counterexample}
Let $H>\sqrt{15}$, $r_H=\sqrt{H^2+1}$, and
\[
v_1=(-1,0),\quad v_2=(1,0),\quad v_3=(0,H).
\]
Set $P_H=\frac13\sum_{k=1}^3\delta_{v_k}$ and $C_H=\conv\{v_1,v_2,v_3\}$.
For $m=2$,
\begin{equation}
V_2(P_H;C_H)=\frac23+\frac{r_H}{2}.
\label{eq:constrained-min}
\end{equation}
The exterior pair $\bm x^\circ=(v_3,(0,-1/\sqrt{15}))$ satisfies
\begin{equation}
V_2(P_H;C_H)-J_{P_H}(\bm x^\circ)
=\gamma_H=\frac{H+4-r_H-\sqrt{15}}6>0.
\label{eq:gap}
\end{equation}
Hence $\Gamma_2(P_H,C_H)\ge\gamma_H$ and every global minimizer has at least one point outside $C_H$.

The key is to compare one exterior pair with all pairs inside the triangle.
To find the constrained minimum, we first show that, for $x,z\in C_H$,
\begin{equation}
\sum_{k=1}^3\bigl(\norm{x-v_k}+\norm{z-v_k}\bigr)
\ge 2+2r_H+\norm{x-z}.
\label{eq:triangle}
\end{equation}
Fix a unit vector $u$ and order the projected vertices as $a\le b\le c$.
For $s=u^\top x$ and $t=u^\top z$ in $[a,c]$, we have $\abs{s-a}+\abs{s-b}+\abs{s-c}=(c-a)+\abs{s-b}$.
Adding the identity for $t$ and using
$\abs{s-b}+\abs{t-b}\ge\abs{s-t}$ gives the projected inequality.
The projected pairwise vertex distances sum to $2(c-a)$.
Integrating over $u=(\cos\theta,\sin\theta)$ and using $\norm w=\frac14\int_0^{2\pi}\abs{u^\top w}\,d\theta$ gives \eqref{eq:triangle}, with constant $\sum_{k<\ell}\norm{v_k-v_\ell}=2+2r_H$.
For $m=2$, \eqref{eq:triangle} and $\diam(C_H):=\sup_{x,z\in C_H}\norm{x-z}=r_H$ imply
\[
J_{P_H}(x,z)\ge\frac{2+2r_H}{3}-\frac{\norm{x-z}}6
\ge\frac23+\frac{r_H}2.
\]
Equality at $(v_3,v_1)$ and $(v_3,v_2)$ proves \eqref{eq:constrained-min}.

For an exterior pair with one point at $v_3$ and the other at $(0,-t)$, $t\ge0$, direct substitution gives
\[
J_{P_H}\bigl(v_3,(0,-t)\bigr)
=\frac{2r_H}{3}-\frac H6+\frac23\sqrt{1+t^2}-\frac t6.
\]
The derivative in $t$ is $2t/(3\sqrt{1+t^2})-1/6$, so the minimum along this ray occurs at $t=1/\sqrt{15}$.
Thus
$J_{P_H}(\bm x^\circ)=2r_H/3-H/6+\sqrt{15}/6$, which gives \eqref{eq:gap}.
The quantity $H+4-\sqrt{H^2+1}$ is strictly increasing and equals $\sqrt{15}$ at $H=\sqrt{15}$, so $\gamma_H>0$.
Existence follows from Lemma~\ref{lem:stability} below, and the strict gap excludes every global minimizer from $C_H^2$.
This argument does not require finding the unconstrained minimum.
The condition $H>\sqrt{15}$ is sufficient; necessity is not established.
\end{example}

To extend the discrete counterexample to a density, we need to control how the objective gap changes under perturbations of the distribution.
The following lemma provides this stability bound and the existence of global minimizers used above.

\begin{lemma}[Existence and stability]
\label{lem:stability}
If $P,Q$ are Borel probability measures on $\R^d$ with finite first moments, $J_P$ is continuous and coercive, $\M_m(P)$ is nonempty and compact, and
\begin{equation}
\sup_{\bm x\in(\R^d)^m}\abs{J_P(\bm x)-J_Q(\bm x)}\le2W_1(P,Q).
\label{eq:stability}
\end{equation}
Here $W_1(P,Q)=\inf\E\norm{Y-Z}$, with the infimum over all couplings $Y\sim P$, $Z\sim Q$, is the Wasserstein distance of order one.
Consequently $\abs{V_m(P;D)-V_m(Q;D)}\le2W_1(P,Q)$ for every nonempty $D\subseteq\R^d$.
If a compact set $C$ contains both supports, then
\begin{equation}
\abs{\Gamma_m(P,C)-\Gamma_m(Q,C)}\le4W_1(P,Q).
\label{eq:gap-stability}
\end{equation}
\end{lemma}
\begin{proof}
Write $M_1(P)=\E\norm Y$.
Since $\sum_{i,j}\norm{x_i-x_j}\le2(m-1)\sum_i\norm{x_i}$, the triangle inequality gives
\begin{equation}
J_P(\bm x)\ge\frac{2}{m^2}\sum_i\norm{x_i}-2M_1(P).
\label{eq:coercivity}
\end{equation}
Continuity and coercivity give a nonempty compact minimizer set.
For any coupling $(Y,Z)$ of $P,Q$,
$\abs{\E\norm{x-Y}-\E\norm{x-Z}}\le\E\norm{Y-Z}$ uniformly in $x$.
Summing and taking the infimum over couplings gives \eqref{eq:stability}; taking infima over $D^m$ gives the value bound.
Applying it to the two values in \eqref{eq:price} proves \eqref{eq:gap-stability}.
\end{proof}

The next example shows that failure is not limited to a distribution on finitely many points.

\begin{example}[A density with full convex support]
\label{ex:continuous}
With $H$ and $C_H$ as in Example~\ref{ex:counterexample}, \eqref{eq:gap-stability} gives, for any probability measure $Q$ supported in $C_H$,
\begin{equation}
\Gamma_2(Q,C_H)\ge\gamma_H-4W_1(Q,P_H).
\label{eq:robust-gap}
\end{equation}
Let $U_{k,\eta}$ be uniform on $C_H\cap B(v_k,\eta)$, where $B(v_k,\eta)=\{x\in\R^2:\norm{x-v_k}<\eta\}$, and let $U_{C_H}$ be uniform on $C_H$.
For $\eta>0$ and $0<\varepsilon<1$, define
\begin{equation}
Q_{\varepsilon,\eta}=\frac{1-\varepsilon}{3}\sum_{k=1}^3U_{k,\eta}
+\varepsilon U_{C_H}.
\label{eq:mixture}
\end{equation}
This absolutely continuous law has support exactly $C_H$ and remains close to $P_H$ for small $\varepsilon$ and $\eta$.
Couple each cap to its vertex and the uniform component to an independent uniform choice of vertex.
The latter distances are at most $r_H$, so
\[
W_1(Q_{\varepsilon,\eta},P_H)\le
b_{\varepsilon,\eta}:=(1-\varepsilon)\eta+\varepsilon r_H.
\]
Consequently, if
\begin{equation}
b_{\varepsilon,\eta}<\gamma_H/4,
\label{eq:smoothing-bound}
\end{equation}
then every global two-point minimizer has a point outside its convex support.
For $H=8$, $\eta=10^{-3}$ and $\varepsilon=10^{-4}$, $\gamma_8\approx0.010793$ and $b_{\varepsilon,\eta}\approx0.001806$ give a gap greater than $0.0035$.
This does not settle containment for the uniform law on a triangle.
\end{example}

Symmetric positive-definite (SPD) matrices, used in covariance analysis and diffusion tensor imaging \citep{DrydenKoloydenkoZhou2009}, give another convex feasibility requirement.
Their sample hull preserves positive definiteness; leaving it need not violate that property.
The next example does, under the Frobenius norm; Riemannian and log-Euclidean distances differ \citep{PennecFillardAyache2006,ArsignyEtAl2006,Schwartzman2016}.

\begin{example}[Positive-definite matrix data]
\label{ex:spd}
We can turn the planar construction into three SPD $2\times2$ observations with an indefinite matrix in a globally optimal pair under the Frobenius norm.

Let $\mathbb S^2$ be the space of symmetric $2\times2$ matrices with the Frobenius inner product $\langle A,B\rangle_F=\operatorname{tr}(A^\top B)$ and norm $\norm A_F=\sqrt{\operatorname{tr}(A^\top A)}$, where $\operatorname{tr}$ denotes the trace.
First embed the planar vertices as $A_k=\operatorname{diag}(v_k)$.
These matrices are not SPD; we will transform them below.
This embedding preserves the hull minimum and the competitor value.
Allowing non-diagonal matrices can only lower the unconstrained minimum, so the strict comparison with the data-hull minimum still holds.
Thus every global minimizer over $(\mathbb S^2)^2$ has a matrix outside the data hull.
Existence follows from Lemma~\ref{lem:stability}.

Choose such an outside matrix $Z$.
Strict separation gives $U\in\mathbb S^2$ with $\norm U_F=1$ and $s\in\R$ such that
\[
\operatorname{tr}(UZ)<s<\min_k\operatorname{tr}(UA_k).
\]
Complete $U$ to an orthonormal basis to obtain an orthogonal linear map $T$ on $\mathbb S^2$ with $[T(A)]_{11}=\operatorname{tr}(UA)$.
Set $\Phi(A)=T(A)+\operatorname{diag}(-s,t)$.
Choose $t$ sufficiently large that $[T(Z)]_{22}+t>0$ and all three $\Phi(A_k)$ are SPD; their first diagonal entries are positive, and their determinants grow with $t$.
Then $\Phi(Z)$ has diagonal entries of opposite signs, so it is indefinite.
The map $\Phi$ is bijective and preserves all Frobenius distances, hence also global optimality.
Thus the transformed data are SPD, but one globally optimal pair contains an indefinite matrix.
We do not claim that every optimal pair for these data does so.

For a numerical illustration, generate $A=LL^\top$ with
$L=\bigl(\begin{smallmatrix}e^\xi&0\\\omega&e^\zeta\end{smallmatrix}\bigr)$,
where $\xi\sim N(0,1.5^2)$, $\zeta\sim N(0,0.25^2)$ and $\omega\sim N(0,1)$ are independent.
Using NumPy's \texttt{default\_rng} with seed 9, we drew all 80 values of $\xi$, then $\zeta$, then $\omega$, giving $N=80$ SPD matrices with smallest input eigenvalue $0.0019016$.
For $m=2$, we applied L-BFGS-B without bounds to the unsmoothed empirical objective over $(\mathbb S^2)^2$, starting from 32 pairs of distinct observed SPD matrices selected using a separate \texttt{default\_rng} with seed 20261014.
Every run returned an indefinite representative; the smallest eigenvalue in the lowest-objective pair was $-0.0026581$, with gradient norms below $1.2\times10^{-7}$ across runs.
This selected numerical example neither establishes global optimality nor estimates a failure frequency.
\end{example}

\section{Empirical consequences and constraint cost}
We next study whether more observations resolve the containment problem and how much the sample-hull constraint increases the minimum objective.
Let $Y_1,Y_2,\ldots$ be i.i.d. from a compactly supported $P$, and define
\[
P_N=\frac1N\sum_{k=1}^N\delta_{Y_k},\qquad
C=\conv(\supp P),\qquad C_N=\conv\{Y_1,\ldots,Y_N\}.
\]
Almost surely, all observations lie in $\supp P$, so $C_N\subseteq C$ for every $N$.
Write $w_N=W_1(P_N,P)$ and $h_N=d_H(C_N,C)$, where the Hausdorff distance is $d_H(A,B)=\max\{\sup_{a\in A}\inf_{b\in B}\norm{a-b},\sup_{b\in B}\inf_{a\in A}\norm{b-a}\}$ for nonempty compact sets $A,B$, and keep $m,d$ fixed.
Almost surely, $w_N\to0$ by empirical Wasserstein convergence \citep{PanaretosZemel2019}, and $h_N\to0$ by density of the observations in the support and a finite covering.
Convexity of $\supp P$ is not required; see \citet{Brunel2019} for random-hull approximation bounds.

\begin{samepage}
\begin{theorem}[Exclusion from the sample hull]
\label{thm:empirical}
Suppose $\Gamma:=\Gamma_m(P,C)>0$ and an empirical configuration $\widehat{\bm x}_N=(\widehat x_{N,1},\ldots,\widehat x_{N,m})\in(\R^d)^m$ satisfies
\begin{equation}
J_{P_N}(\widehat{\bm x}_N)\le V_m(P_N)+a_N,\qquad a_N\ge0.
\label{eq:approx}
\end{equation}
Then
\begin{equation}
4w_N+a_N<\Gamma\quad\Longrightarrow\quad
\widehat{\bm x}_N\notin C_N^m.
\label{eq:exclusion}
\end{equation}
If $a_N\to0$ almost surely, every configuration satisfying \eqref{eq:approx} has at least one point outside $C_N$ eventually almost surely.
\end{theorem}
\end{samepage}
\begin{proof}
Lemma~\ref{lem:stability} and \eqref{eq:approx} give
\begin{equation}
J_P(\widehat{\bm x}_N)-V_m(P)\le4w_N+a_N.
\label{eq:excess}
\end{equation}
If $\widehat{\bm x}_N\in C_N^m\subseteq C^m$, the left-hand side is at least $\Gamma$.
This proves \eqref{eq:exclusion}.
Almost-sure convergence of $w_N$ and $a_N$ gives the final assertion, simultaneously for all configurations with the specified error bound.
\end{proof}

Theorem~\ref{thm:empirical} shows that more observations do not restore containment when the population gap is positive and the global objective error tends to zero.
Here $a_N$ is zero for exact minimizers; convergence to a local solution alone does not provide this bound.
For measurable configurations and error bounds,
$\Prb\{\widehat{\bm x}_N\in C_N^m\}\le\Prb\{4w_N+a_N\ge\Gamma\}$,
which can be controlled using Wasserstein concentration bounds \citep{FournierGuillin2015}.
Example~\ref{ex:continuous} gives $\Gamma\ge\gamma_H-4b_{\varepsilon,\eta}>0$ under \eqref{eq:smoothing-bound}.

A sample-hull constraint guarantees feasibility in any convex region containing the observations, but it can increase the objective.
The next result relates this increase to the population cost.
Selecting observed representatives also ensures containment, as in support-point resampling \citep{HuangJosephMak2022}.

\begin{proposition}[Constrained-value bound]
\label{prop:price}
Set $\Gamma_N=V_m(P_N;C_N)-V_m(P_N)$ and $\Gamma=\Gamma_m(P,C)$.
Then
\begin{equation}
-4w_N\le\Gamma_N-\Gamma
\le4w_N+2(1-1/m)h_N.
\label{eq:price-bound}
\end{equation}
In particular, $\Gamma_N\to\Gamma$ almost surely.
This limiting cost is zero if and only if some population global minimizer lies in $C^m$.
\end{proposition}
\begin{proof}
Let $A_N=V_m(P_N;C_N)$ and $A=V_m(P;C)$.
Since $C_N\subseteq C$, Lemma~\ref{lem:stability} gives $A_N\ge A-2w_N$.
For the upper bound, take a minimizer $\bm x^*$ over $C^m$ and project each coordinate onto $C_N$, obtaining $z_i$ with $\norm{z_i-x_i^*}\le h_N$.
Projection cannot increase distance to a sample point, so the attraction term of $J_{P_N}$ cannot increase.
The repulsion term increases by at most
\[
\frac1{m^2}\sum_{i\ne j}
\bigl(\norm{x_i^*-z_i}+\norm{x_j^*-z_j}\bigr)
\le2(1-1/m)h_N.
\]
Thus $A_N\le A+2w_N+2(1-1/m)h_N$.
Subtracting the unconstrained values and using
$\abs{V_m(P_N)-V_m(P)}\le2w_N$ gives the bounds.
The final assertion follows because the constrained minimum is attained on compact $C^m$.
\end{proof}

Sample-hull containment provides a sufficient, potentially conservative guarantee of feasibility under convex constraints.
The unconstrained energy-distance objective does not provide this guarantee, and the examples show that its minimizers can also violate the underlying constraints.
Imposing the hull constraint makes feasibility explicit.
Proposition~\ref{prop:price} quantifies its objective cost, which can remain positive as the sample grows for fixed $m$.

\section*{Funding}
Zhaohui Li's research is supported by NSFC 12401327.

\section*{Declaration of generative AI and AI-assisted technologies in the writing process}
The author used a large language model to assist with language and writing review, the presentation and formatting of figures, and the search, collection, and organization of references.
The author independently verified the manuscript, figures, references, and citations, reviewed and edited all resulting material, and takes full responsibility for the content of the publication.

\setlength{\bibsep}{0pt}
\end{document}